\documentclass[review]{elsarticle}
\usepackage{latexsym,amsfonts,amsmath,epsfig,amssymb,mathrsfs,enumerate}
\usepackage[pdftex, bookmarksnumbered, bookmarksopen, colorlinks, citecolor=blue, linkcolor=blue]{hyperref}
\usepackage{algorithm} 
\usepackage{algorithmic} 
\usepackage{booktabs}
\usepackage{lineno,hyperref}

\journal{Engineering Optimization}

\begin{document}

\begin{frontmatter}

\title{Iterative Complexity  and Applications of Linearized Generalized Alternating Direction Method of Multipliers with Multi-block Case}

\author{Jian He $^{a, b}$}
\author{Bangzhong Zhang $^{c}$}
\author{Jinlin Li $^{a^{*}, d}$}
\cortext[mycorrespondingauthor]{Corresponding author}
\ead{20102197@cqu.edu.cn}

\address[mainaddress]{School of Economics and Business Administration, Chongqing University, Chongqing 400030, People's Republic of China.}
\address[mysecondaryaddress]{Chongqing Chemical Industry Vocational College, Chongqing 401228, People's Republic of China.}
\address[mythirdaryaddress]{Chongqing Tsinghua High School, Chongqing 400054, People's Republic of China.}
\address[myforthdaryaddress]{Postdoctoral Workstation of Chongqing Rural Commercial Bank, Chongqing 400023, People's Republic of China.}

\begin{abstract}
We consider a multi-block separable convex optimization problem with the linear constraints, where the objective function is the sum of $m$ individual convex functions without overlapping variables. The linearized version of the generalized alternating direction method of multipliers (L-GADMM) is particularly efficient for the two-block separable convex programming problem and its convergence was proved when two blocks of variables
are alternatively updated. However, the convergence and some practical applications of the extension ($m\geq3$) of the L-GADMM is still in its infancy. In this paper, we extend this algorithm to the general case where the objective function consists of the sum of $m$-block convex functions. Theoretically, we prove global convergence of the new method and establish the worst-case convergence rate in the ergodic and nonergodic senses for the proposed algorithm. The efficiency of the new method is further demonstrated through numerical results on the calibration of the correlation matrices.
\end{abstract}


\end{frontmatter}


\section{Introduction}
We consider the multiple-block separable convex minmization problem with linear constraints.
\begin{equation}
\mathrm{min} \left\{\sum_{i=1}^{m}\vartheta_{i}(x_{i})|\sum_{i=1}^{m}\mathcal{A}_{i}x_{i}=b, x_{i}\in\mathcal{X}_{i},i=1, 2, ..., m\right\},
\end{equation}
where $m$ is any positive integer, $\vartheta_{i} : \mathcal{R}^{n_{i}}\rightarrow\mathcal{R} (i=1, 2, \cdots, m)$ is convex function but not necessarily smooth. Let $\ell$ be also a positive integer, $\mathcal{A}_{i}\in\mathcal{R}^{\ell\times n_{i}} (i=1, 2, \cdots, m)$ be given matrices, $b\in\mathcal{R}^{\ell}$ be a given vector, and $\mathcal{X}_{i}\subseteq\mathcal{R}^{n_{i}} (i=1, 2, \cdots ,m)$ be nonempty closed convex set. Throughout this paper, we let $n=\sum_{i=1}^{m}n_{i}$, $\textbf{x}:=(x_{1},x_{2},\cdots,x_{m})$ and $\mathcal{X}:=\prod_{i=1}^{m}\mathcal{X}_{i}$ and the matrice $\mathcal{A}_{i} (i=1, 2, \cdots , m)$ is assumed to be full column rank. Moreover, the solution set of (1) is assumed to be nonempty. The problem (1) has many numerous applications, e.g., the robust principal component analysis (PCA) model with noisy and incomplete data[1], signal or image processing [2, 3, 4], the quadratic discriminant analysis model [5], statistical learning [1, 6, 7, 8], network optimization [9, 10, 11], DSL dynamic spectrum management problem [12, 13].\\
\indent Though the model (1) has a simple form, the efficient solutions is usually very challenging in practice because of the nonsmoothness of the objective functions and the process of high dimensionality arising from some real-world applications. Thus, it is efficient by the structure-utilizing algorithms for (1) to overcome the difficulties caused by the aforementioned characteristics. Based on the special structure, some operator splitting methods have been developed for the problem (1) with m=2. These methods involve the Peaceman-Rachford splitting method (PRSM) [5, 16], the alternating direction method of multipliers (ADMM) [14, 15], and their variants [17-26]. In fact, to solve the problem (1) with $m=2$, a benchmark is the ADMM proposed originally in [14, 15]. And the schemes of the ADMM for the optimization problem (1) with two-block are written as follows
\begin{equation}
\begin{split}
 \left\{
\begin{aligned}
x_{1}^{k+1} &= \mathrm{argmin}\{\vartheta_{1}(x_{1})-x_{1}^{\top}\mathcal{A}_{1}^{\top}y^{k}+\frac{\rho}{2}\| \mathcal{A}_{1}x_{1}+\mathcal{A}_{2}x_{2}^{k}-b \|^{2}|x_{1}\in\mathcal{X}_{1}\} , \\
x_{2}^{k+1} &= \mathrm{argmin}\{\vartheta_{2}(x_{2})-x_{2}^{\top}\mathcal{A}_{2}^{\top}y^{k}+\frac{\rho}{2}\| \mathcal{A}_{1}x_{1}^{k+1}+\mathcal{A}_{2}x_{2}-b \|^{2}|x_{2}\in\mathcal{X}_{2}\} , \\
y^{k+1}&=y^{k}-\rho(\mathcal{A}_{1}x_{1}^{k+1}+\mathcal{A}_{2}x_{2}^{k+1}-b),\\
\end{aligned}
\right.
\end{split}
\end{equation}
where $y\in\mathcal{R}^{\ell}$ is the lagrange multiplier and $\rho>0$ is a penalty parameter. Gabay [27] pointed out that the ADMM was an application of the well-known Douglas-Rachford splitting method (DRSM) in [16] to the dual of (2); and as an
application of the proximal point algorithm (PPA) in [28], Eckstein and Bertsekas [23]  proposed the following
generalized ADMM (GADMM) scheme:
\begin{equation}
\begin{split}
\left\{
\begin{aligned}
x_{1}^{k+1} &= \mathrm{argmin}\{\vartheta_{1}(x_{1})-x_{1}^{\top}\mathcal{A}_{1}^{\top}y^{k}+\frac{\rho}{2}\| \mathcal{A}_{1}x_{1}+\mathcal{A}_{2}x_{2}^{k}-b \|^{2}|x_{1}\in\mathcal{X}_{1}\} , \\
x_{2}^{k+1} &=\mathrm{argmin}\{\vartheta_{2}(x_{2})-x_{2}^{\top}\mathcal{A}_{2}^{\top}y^{k}+\frac{\rho}{2}\| \gamma \mathcal{A}_{1}x_{1}^{k+1}+(1-\gamma)(b-\mathcal{A}_{2}x_{2}^{k})\\
&~~+\mathcal{A}_{2}x_{2}-b \|^{2}|x_{2}\in\mathcal{X}_{2}\} , \\
y^{k+1}&=y^{k}-\rho\left( \gamma \mathcal{A}_{1}x_{1}^{k+1}+(1-\gamma)(b-\mathcal{A}_{2}x_{2}^{k})+\mathcal{A}_{2}x_{2}^{k+1}-b\right),\\
\end{aligned}
\right.
\end{split}
\end{equation}
where the parameter $\gamma\in(0,2)$ is a relaxation factor. If the parameter $\gamma=1$, the generalized scheme (3) will reduce to the original ADMM scheme (2). Based on the above iterative scheme (3), {Fang} et al. [29] recently proposed the following linearized versions of GADMM (``L-GADMM'' for short) for (3) by adding a proximal term to its subproblem:
\begin{equation}
\begin{split}
\left\{
\begin{aligned}
x_{1}^{k+1}=& \mathrm{argmin}\{\vartheta_{1}(x_{1})-x_{1}^{\top}\mathcal{A}_{1}^{\top}y^{k}+\frac{\rho}{2}\| \mathcal{A}_{1}x_{1}+\mathcal{A}_{2}x_{2}^{k}-b \|^{2}\\
&+\frac{1}{2}\|x_{1}-x_{1}^{k}\|_{P_{1}}^{2}|x_{1}\in\mathcal{X}_{1}\} , \\
x_{2}^{k+1}=&\mathrm{argmin}\{\vartheta_{2}(x_{2})-x_{2}^{\top}\mathcal{A}_{2}^{\top}y^{k}+\frac{\rho}{2}\|
\gamma \mathcal{A}_{1}x_{1}^{k+1}+(1-\gamma)(b-\mathcal{A}_{2}x_{2}^{k})\\
&+\mathcal{A}_{2}x_{2}-b \|^{2}|x_{2}\in\mathcal{X}_{2}\} , \\
y^{k+1}=&y^{k}-\rho( \gamma \mathcal{A}_{1}x_{1}^{k+1}+(1-\gamma)(b-\mathcal{A}_{2}x_{2}^{k})+\mathcal{A}_{2}x_{2}^{k+1}-b),\\
\end{aligned}
\right.
\end{split}
\end{equation}
and
\begin{equation}
\begin{split}
\left\{
\begin{aligned}
x_{1}^{k+1}=&\mathrm{argmin}\{\vartheta_{1}(x_{1})-x_{1}^{\top}\mathcal{A}_{1}^{\top}y^{k}+\frac{\rho}{2}\| \mathcal{A}_{1}x_{1}+\mathcal{A}_{2}x_{2}^{k}-b \|^{2}\\
&+\frac{1}{2}\|x_{1}-x_{1}^{k}\|_{P_{1}}^{2}|x_{1}\in\mathcal{X}_{1}\} , \\
x_{2}^{k+1}=&\mathrm{argmin}\{\vartheta_{2}(x_{2})-x_{2}^{\top}\mathcal{A}_{2}^{\top}y^{k}+\frac{\rho}{2}\| \gamma \mathcal{A}_{1}x_{1}^{k+1}+(1-\gamma)(b-\mathcal{A}_{2}x_{2}^{k})\\
&+\mathcal{A}_{2}x_{2}-b \|^{2}+\frac{1}{2}\|x_{2}-x_{2}^{k}\|_{P_{2}}^{2}|x_{2}\in\mathcal{X}_{2}\} , \\
y^{k+1}=&y^{k}-\rho( \gamma \mathcal{A}_{1}x_{1}^{k+1}+(1-\gamma)(b-\mathcal{A}_{2}x_{2}^{k})+\mathcal{A}_{2}x_{2}^{k+1}-b),\\
\end{aligned}
\right.
\end{split}
\end{equation}
where $P_{1}\in \Re^{n_{1}\times n_{1}}$ and $P_{2}\in \Re^{n_{2}\times n_{2}}$ are symmetric positive definite matrices. In this paper, we focus on the above linearized version of the GADMM, which is an efficient and simple acceleration algorithm. Note that we use the notation $\|x_{1}\|_{P_{1}}$ to denote the quantity $\sqrt{x_{1}^{\top}P_{1}x_{1}}$. If $\mathcal{X}_{1}=\Re^{n_{1}}$ and $P_{1}=\tau I_{n_{1}}-\rho \mathcal{A}_{1}^{\top}\mathcal{A}_{1}$ with the requirement $\tau>\rho \|\mathcal{A}_{1}^{\top}\mathcal{A}_{1}\|_{2}$, where $\|\cdot\|_{2}$ denotes the spectral norm of a matrix, then the $x_{1}$-subproblem in (4) or (5) reduces to estimating the resolvent operator of $ \vartheta_{1}$:
\begin{equation*}
x_{1}^{k+1}=(I+\frac{1}{\tau}\partial \vartheta_{1})^{-1}(\textbf{t})=\mathrm{argmin}\{\vartheta_{1}(x_{1})+\frac{\tau}{2}\|x_{1}-\textbf{t}\|^{2}\},
\end{equation*}
where $\partial(\cdot)$ denotes the subdifferential of a convex function and $\textbf{t}=\frac{1}{\tau}((\tau I_{n_{1}}-\rho \mathcal{A}_{1}^{\top}\mathcal{A}_{1})x_{1}^{k}-\rho \mathcal{A}_{1}^{\top}\mathcal{A}_{2}x_{2}^{k}+\mathcal{A}_{1}^{\top}y^{k}+\rho \mathcal{A}_{1}^{\top}b))$. Therefore, the scheme (4) includes the linearized version of the ADMM (see [26, 30, 31]) as a special case with $P_{1}=\tau I_{n_{1}}-\rho \mathcal{A}_{1}^{\top}\mathcal{A}_{1}$ and $\gamma=1$.\\
\indent Theoretically, {Fang} et al. [29] showed the worst $\mathcal{O}(1/t)$ convergence rate in both ergodic and nonergodic senses for the L-GADMM (4) and (5). Numerically, {Fang} et al. [29] verified that both the L-GADMM (4) and the scheme (5) were quite efficient for solving some rather new and core applications arising in statistical learning. However, both the iterative scheme (4) and (5) are only suit for the minimization  problem (1) with $m=2$. In fact, there are some efficient algorithms to solve the problem (1) with $m\geq3$, we thus refer to [20, 32, 33].  In this paper, we regroup the $m$ variables $x_{i}$ $(i=1, 2, \cdots , m)$ as two blocks $(x_{1}, x_{2}, \cdots , x_{m-1})$ and $x_{m}$ and further separate the subproblems related to the first block $(x_{1}, x_{2}, \cdots , x_{m-1})$ into $m-1$ parallel smaller subproblems. Furthermore, our algorithm is written as follows
\begin{equation}
\begin{split}
\left\{
\begin{aligned}
x_{j}^{k+1} &= \mathrm{argmin}\{\vartheta_{j}(x_{j})+\frac{\rho}{2}\|\mathcal{A}_{j}x_{j}+\sum_{i=1,i\neq j}^{m}\mathcal{A}_{i}x_{i}^{k}-b-\frac{y^{k}}{\rho}\|^{2}+\frac{1}{2}\|x_{j}-x_{j}^{k}\|_{P_{j}}^{2}|x_{j}\in\mathcal{X}_{j}\} , \\
&~~~~~~~~~~~~~~~~~~~~~~~~~~~~~~~~j=1, ~\cdots~ ,~ m-1,\\
x_{m}^{k+1} &=\mathrm{argmin}\{\vartheta_{m}(x_{m})+\frac{\rho}{2}\| \gamma \sum_{i=1}^{m-1}\mathcal{A}_{i}x_{i}^{k+1}+(1-\gamma)(b-\mathcal{A}_{m}x_{m}^{k})+\mathcal{A}_{m}x_{m}-b-\frac{y^{k}}{\rho} \|^{2}\\
&~~+\frac{1}{2}\|x_{m}-x_{m}^{k}\|_{P_{m}}^{2}|x_{m}\in\mathcal{X}_{m}\} , \\
y^{k+1}&=y^{k}-\rho\left( \gamma \sum_{i=1}^{m-1}\mathcal{A}_{i}x_{i}^{k+1}+(1-\gamma)(b-\mathcal{A}_{m}x_{m}^{k})+\mathcal{A}_{m}x_{m}^{k+1}-b)\right),\\
\end{aligned}
\right.
\end{split}
\end{equation}
where $P_{i}\in\mathcal{R}^{n_{i}\times n_{i}} (i=1, 2, \cdots, m-1)$ are symmetric positive definite matrices, $\gamma\in(0,2)$ is a relaxation factor and $P_{m}\in\mathcal{R}^{n_{m}\times n_{m}}$ is a symmetric matrix with $P_{m}+ \frac{\rho}{\gamma}\mathcal{A}_{m}^{\top}\mathcal{A}_{m}\succ\textbf{0}$.  Here for a matrix $P\succ\textbf{0}$ means that $P$ is positive definite. It is note worthy that the above iterative scheme (6) has the following features: (I) It minimizes $m$ subproblems with respect to $x_{i} (i=1, 2, \cdots, m)$ separately and solves the first $m-1$ subproblems in parallel; (II) The variable $x_{m}$ is updated by using all the latest information $(x_{1}^{k+1}, x_{2}^{k+1}, \cdots, x_{m-1}^{k+1})$. In fact, the scheme (6) is a partially parallel splitting method. (III) All the subproblems for the iterative scheme (6) are linearized and the corresponding matrices for the quadratic proximal terms are positive definite but the matrix $P_{m}$. Specially, the matrix $P_{m}$ is not necessarily positive definite, we only need to ensure that the matrix satisfies the condition:  $P_{m}+\frac{\rho}{\gamma}\mathcal{A}_{m}^{\top}\mathcal{A}_{m}\succ\textbf{0}$.\\
\indent \textbf{Our contribution:}\\
\indent $\bullet$ We extend L-GADMM (6) to the general case where the objective function consists of the sum of $m$-block ($m\geq 2$) convex functions. \\
\indent $\bullet$ We show the global convergence of the new method and establish the worst-case convergence rate in the ergodic and nonergodic senses for the proposed algorithm (6).  \\
\indent $\bullet$  The efficiency of the new method is further demonstrated through the calibration of the correlation matrices.\\
\indent The rest of this paper is organized as follows. In section 2, we give some notations and characterize problem (1) by a variational inequality problem. Then, we propose L-GADMM iterative scheme and prove its global convergence in section 3. In section 4, we establish a worst-case $\mathcal{O}(1/t)$ convergence rate in both ergodic and nonergodic senses for the proposed method. In Section 5, we apply the new algorithm to solve the matrix optimization problem and report some experiment results. Finally, we draw some conclusions in Section 6.

\section{Preliminaries}
In this section, we give some notations and summarize some properties that will be used throughout this paper.\\
\indent We use $\mathcal{R}^{n}$ to denote an $n$-dimensional Euclidean space. The superscript ``$^{\top}$" represents the transpose operation for matrix/vector variables. For any a vector $x\in\mathcal{R}^{n}$, let $\|x\|_{p}=\sqrt{x^{\top}x}$ be the $p$-norm. Moreover, $I$ represents the identity matrix and $\textbf{0}$ represents a zero vector or matrix. It is well known that solving (1) is equivalent to solving a variational inequality problem (denoted by $VI(\mathcal{W},\mathcal{F},\vartheta)$): Find a vector $w^{*}\in \mathcal{W}$ such that
\begin{equation}
\vartheta(u)-\vartheta(u^{*})+(w-w^{*})^{\top}\mathcal{F}(w^{*})\geq0,~~~~\forall w\in \mathcal{W},
\end{equation}
where
\begin{equation}       
u=\left (                 
  \begin{array}{ccc}   
    x_{1}\\  
    x_{2}\\
    \vdots\\
    x_{m}\\
  \end{array}\right),~~
w=\left (                 
  \begin{array}{ccc}   
    x_{1}\\  
    x_{2}\\
    \vdots\\
    x_{m}\\
    y\\  
  \end{array}\right)
,~~~~\vartheta(u)=\sum_{i=1}^{m}\vartheta_{i}(x_{i})
  \end{equation}
  and
  \begin{equation}
\mathcal{F}(w)=\left (                 
  \begin{array}{ccc}   
    -\mathcal{A}_{1}^{\top}y\\  
    -\mathcal{A}_{2}^{\top}y\\  
    \vdots\\
     -\mathcal{A}_{m}^{\top}y\\
     \sum_{i=1}^{m}\mathcal{A}_{i}x_{i}-b
  \end{array}
\right)
,~~ ~~\mathcal{W}=\mathcal{X}_{1}\times\mathcal{X}_{2}\times\cdots\times\mathcal{X}_{m}\times\Re^{\ell},
\end{equation}
The solution set of $VI(\mathcal{W},\mathcal{F},\vartheta)$ denoted $\mathcal{W}^{*}$. Obviously, $\mathcal{W}^{*}$ is nonempty since the solution set of (1) is nonempty. It is easy to verify that the mapping $\mathcal{F}(\cdot)$ defined by (9) is monotone and satisfies the following a property:
\begin{equation*}
(w_{1}-w_{2})^{\top}(\mathcal{F}(w_{1})-\mathcal{F}(w_{2})\geq0,~~~~\forall ~~w_{1}, w_{2}\in\mathcal{W}.
\end{equation*}
\indent In the following, we will give some matrices which will be used later. Set
\begin{equation}
\mathcal{Q}=\left(                 
  \begin{array}{ccc}   
    \mathcal{G}_{1} & \textbf{0} & \textbf{0}\\  
    \textbf{0} & \rho \mathcal{A}_{m}^{\top}\mathcal{A}_{m}+P_{m} & (1-\gamma)\mathcal{A}_{m}^{\top}\\  
    \textbf{0}&-\mathcal{A}_{m}&\frac{1}{\rho} I_{\ell}\\
  \end{array}
\right),
\end{equation}
\begin{equation}
\mathcal{M}=\left(                 
  \begin{array}{ccc}   
    I_{1} & \textbf{0} & \textbf{0}\\  
    \textbf{0} & I_{n_{m}} & \textbf{0}\\  
    \textbf{0}&-\rho \mathcal{A}_{m}&\gamma I_{\ell}\\
  \end{array}
\right),
\end{equation}
and
\begin{equation}
\mathcal{H}=\left(                 
  \begin{array}{ccc}   
    \mathcal{G}_{1} & \textbf{0} & \textbf{0}\\  
    \textbf{0} & P_{m}+\frac{\rho}{\gamma}\mathcal{A}_{m}^{\top}\mathcal{A}_{m}& \frac{1-\gamma}{\gamma}\mathcal{A}_{m}^{\top}\\  
    \textbf{0}&\frac{1-\gamma}{\gamma}\mathcal{A}_{m}&\frac{1}{\gamma\rho} I_{\ell}\\
  \end{array}
\right),
\end{equation}
where
\begin{equation*}
I_{1}=\left(                 
  \begin{array}{ccccc}   
I_{n_{1}}&\textbf{0}&\cdots&\textbf{0}\\
\textbf{0}&I_{n_{2}}&\cdots&\textbf{0}\\
\cdots&\cdots&\cdots&\cdots\\
\textbf{0}&\textbf{0}&\cdots&I_{n_{m-1}}\\
  \end{array}
\right),~~
\mathcal{G}_{1}=\left(                 
  \begin{array}{ccccc}   
    P_{1} & -\rho \mathcal{A}_{1}^{\top}\mathcal{A}_{2} &\cdots& -\rho\mathcal{A}_{1}^{\top}\mathcal{A}_{m-1}\\  
    -\rho\mathcal{A}_{2}^{\top}\mathcal{A}_{1} & P_{2} &\cdots& -\rho\mathcal{A}_{2}^{\top}\mathcal{A}_{m-1}\\
    \cdots&\cdots&\cdots&\cdots\\
    -\rho\mathcal{A}_{m-1}^{\top}\mathcal{A}_{1}&-\rho\mathcal{A}_{m-1}^{\top}\mathcal{A}_{2}&\cdots&P_{m-1}\\
  \end{array}
\right).
\end{equation*}
\indent Note that the matrix $\mathcal{G}_{1}$ is symmetric but not necessarily positive definite. In fact, it is not difficult to satisfy the positive definiteness of the matrix  $\mathcal{G}_{1}$. For instance, for the matrix
$\mathcal{G}_{1}=\left(
  \begin{array}{ccccc}
P_{1}&-\rho \mathcal{A}_{1}^{\top}\mathcal{A}_{2}\\
-\rho\mathcal{A}_{2}^{\top}\mathcal{A}_{1}&P_{2}
  \end{array}
  \right)$,
 then the $\mathcal{G}_{1}$ is positive definite when both $P_{1}-\rho\|A_{1}\|^{2}\cdot I_{n_{1}}$ and $P_{2}-\rho\|A_{2}\|^{2}\cdot I_{n_{2}}$ are positive definite. \\
 \indent Then, we will give the following a lemma to elaborate some properties of the three matrices $\mathcal{M}$, $\mathcal{Q}$ and $\mathcal{H}$. It is easy to obtain the following results by the matrices $\mathcal{M}$, $\mathcal{Q}$, $\mathcal{H}$ just defined, we thus omit their proofs.  \\
 \textbf{Lemma 2.1} Let the matrices $\mathcal{Q}, \mathcal{M}, \mathcal{H}$ be defined in (10), (11) and (12), respectively. If $\alpha\in (0,2)$ and $\mathcal{G}_{1}\succ\textbf{0}$, then \\
~~~~(I).
$\mathcal{Q}=\mathcal{H}\mathcal{M}$.\\
~~~~(II).  The matrix $\mathcal{H}$ is positive definite. \\                
~~~~(III). $\mathcal{N}:=\mathcal{Q}^{\top}+\mathcal{Q}-\mathcal{M}^{\top}\mathcal{H}\mathcal{M}\succ\textbf{0}$.\\
\section{The global convergence of the L-GADMM  scheme}
In this section, we further develop the iterative scheme (6) for solving (1) and establish its global convergence.
\begin{table}[H]
\begin{center}
 \begin{tabular}{lcl}
  \toprule
  $\textbf{Algorithm 1}$ A linearized GADMM for $VI (W, \mathcal{F}, \vartheta)$~~~~~~~~~~\\
  \midrule
  Input $\gamma\in(0,2)$, $\rho>0$, $P_{i}\in\mathcal{R}^{n_{i}\times n_{i}} (i=1,2,\cdots,m-1)$ such that the matrix $\mathcal{G}_{1}$ is positive\\ definite and $P_{m}\in\mathcal{R}^{n_{m}\times n_{m}}$ is a symmetric matrix with $P_{m}+\frac{\rho}{\gamma}\mathcal{A}_{m}^{\top}\mathcal{A}_{m}\succ\textbf{0}$.  \\
  Initialize $(x_{1},x_{2},\cdots,x_{m},y)=(x_{1}^{0},x_{2}^{0},\cdots,x_{m}^{0},y^{0}), k=0$.\\
  $\textbf{while}$ a stopping criterion is not satisfied \textbf{do}\\
  ~~(1)Compute $w^{k+1}=(x_{1}^{k+1},x_{2}^{k+1},\cdots,x_{m}^{k+1},y^{k+1})$ by (6).\\
  ~~(2)Set $k=k+1$.\\
   $\textbf{end while}$ \\
   Output $x_{1}^{k+1},x_{2}^{k+1},\cdots,x_{m}^{k+1}$.\\
  \bottomrule
 \end{tabular}
\end{center}
\end{table}

\indent Throughtout this paper, we definite the following an auxiliary sequence $\{\bar{w}^{k}\}$, which is used to prove the global convergence of Algorithm 1.
\begin{equation}       
\bar{w}^{k}=\left(                 
  \begin{array}{ccc}   
    \bar{x}_{1}^{k}\\  
    \bar{x}_{2}^{k}\\  
    \vdots\\
    \bar{x}_{m}^{k}\\
    \bar{y}^{k}\\
  \end{array}
\right)=\left(
  \begin{array}{ccc}   
    x_{1}^{k+1}\\  
    x_{2}^{k+1}\\   
    \vdots\\
    x_{m}^{k+1}\\
    y^{k}-\rho(\sum_{i=1}^{m-1}\mathcal{A}_{i}x_{i}^{k+1}+\mathcal{A}_{m}x_{m}^{k}-b)\\
  \end{array}
\right),               
\end{equation}
where $(x_{1}^{k+1},x_{2}^{k+1},\cdots,x_{m}^{k+1})$ is generated by the scheme (6) from $(x_{1}^{k},x_{2}^{k},\cdots,x_{m}^{k})$. Furthermore, according to (10) and (13), we get
\begin{equation}
w^{k+1}=w^{k}-\mathcal{M}(w^{k}-\bar{w}^{k}).
\end{equation}
\indent Now, we give the convergence analysis of Algorithm 1. \\
\textbf{Lemma 3.1} Let the sequence $\{w^{k}\}$ be generated by L-GACMM scheme (6). Then we have
\begin{equation}
\vartheta(u)-\vartheta(\bar{u}^{k})+(w-\bar{w}^{k})^{\top}\mathcal{F}(\bar{w}^{k})\geq(w-\bar{w}^{k})^{\top}\mathcal{Q}(w^{k}-\bar{w}^{k}),~~~~\forall w\in \mathcal{W},
\end{equation}
where $\mathcal{F}$, $\mathcal{Q}$ is defined in (9) and (10), respectively.\\
$Proof$ By the first-order optimality condition of the $x_{j}$-subproblem $(j=1,...,m-1)$ in (6), for any $x_{j}\in\mathcal{X}_{j}$, we have
\begin{small}
\begin{equation}
\vartheta_{j}(x_{j})-\vartheta_{j}(x_{j}^{k+1})+(x_{j}-x_{j}^{k+1})^{\top}\left\{-\mathcal{A}_{j}^{\top}\left(y^{k}-\rho(\mathcal{A}_{j}x_{j}^{k+1}+\sum_{i=1,i\neq j}^{m}\mathcal{A}_{i}x_{i}^{k}-b)\right)+P_{j}(x_{j}^{k+1}-x_{j}^{k})\right\}\geq0.
\end{equation}
\end{small}
It follows from (13) and (16) that
 \begin{equation}
 \begin{split}
\vartheta_{j}(x_{j})-\vartheta_{j}(\bar{x}_{j}^{k})+(x_{j}-\bar{x}_{j}^{k})^{\top}\left\{-\mathcal{A}_{j}^{\top}\bar{y}^{k}+\rho\mathcal{A}_{j}^{\top}\left(\sum_{i=1,i\neq j}^{m-1}\mathcal{A}_{i}(x_{i}^{k}-\bar{x}_{i}^{k})\right)+P_{j}(\bar{x}_{j}^{k}-x_{j}^{k})\right\}\\
\geq0, ~~\forall x_{j}\in\mathcal{X}_{j}.
\end{split}
\end{equation}
Furthermore, by the optimality of $x_{m}$-subproblem in (6), for any $x_{m}\in\mathcal{X}_{m}$, we obtain
\begin{small}
\begin{equation}
\begin{split}
&\vartheta_{m}(x_{m})-\vartheta_{m}(x_{m}^{k+1})+\left(x_{m}-x_{m}^{k+1}\right)^{\top}\\
&\times\left\{-\mathcal{A}_{m}^{\top}\left(y^{k}-\rho\left(\gamma\sum_{j=1}^{m-1}\mathcal{A}_{j}x_{j}^{k+1}+(1-\gamma)(b-\mathcal{A}_{m}x_{m}^{k})+\mathcal{A}_{m}x_{m}^{k+1}-b\right)\right)+P_{m}(x_{m}^{k+1}-x_{m}^{k})\right\}\geq0.
\end{split}
\end{equation}
\end{small}
From (13) and (18), we have
\begin{equation}
\begin{split}
\vartheta_{m}(x_{m})-\vartheta_{m}(\bar{x}_{m}^{k})&+(x_{m}-\bar{x}_{m}^{k})^{\top}(-\mathcal{A}_{m}^{\top}\bar{y}^{k}-(1-\gamma)\mathcal{A}_{m}^{\top}(y^{k}-\bar{y}^{k})\\
&+(P_{m}+\rho\mathcal{A}_{m}^{\top}\mathcal{A}_{m})(\bar{x}_{m}^{k}-x_{m}^{k}))\geq0,~~\forall x_{m}\in\mathcal{X}_{m}.
\end{split}
\end{equation}
Moreover, by (13), we get
\begin{equation}
\left(y-\bar{y}^{k}\right)^{\top}\left\{\left(\sum_{j=1}^{m}\mathcal{A}_{j}\bar{x}_{j}^{k}-b\right)-\mathcal{A}_{m}\left(\bar{x}_{m}^{k}-x_{m}^{k}\right)+\frac{1}{\rho}\left(\bar{y}^{k}-y^{k}\right)\right\}\geq0, \forall y\in\mathcal{R}^{\ell}.
\end{equation}
It follows from (8), (9), (10), (17), (19) and (20) that
\begin{equation}
\vartheta(u)-\vartheta(\bar{u}^{k})+\left(w-\bar{w}^{k}\right)^{\top}\left\{\mathcal{F}(\bar{w}^{k})+\mathcal{Q}\left(\bar{w}^{k}-w^{k}\right)\right\}\geq0,~~~~\forall w\in \mathcal{W},
\end{equation}
which implies (15).$~~~~~~~~~~~~~~~~~~~~~~~~~~~~~~~~~~~~~~~~~~~~~~~~~~~~~~~~~~~~~~~~~~~~~~~~~~~~~~~~~~~~~~~~~~~~~~~~\square$\\
\indent From (21), if $w^{k}=\bar{w}^{k}$, we have
\begin{equation}
\vartheta(u)-\vartheta(\bar{u}^{k})+\left(w-\bar{w}^{k}\right)^{\top}\mathcal{F}(\bar{w}^{k})\geq0,~~~~\forall w\in \mathcal{W}.
\end{equation}
It follows from (22) that $\bar{w}^{k}=(\bar{x}_{1}^{k}, \bar{x}_{2}^{k},..., \bar{x}_{m}^{k}, \bar{y}^{k})$ is a solution of $VI(\mathcal{W},\mathcal{F},\vartheta)$. Moreover, by (13), we can obtain the following a result which implies the Algorithm 1 to produce two infinite sequences $\{w^{k}\}$ and $\{\bar{w}^{k}\}$:
\\
\textbf{$Remark~ 3.1$} If ${A}_{i}x_{i}^{k}={A}_{i}\bar{x}_{i}^{k}$ $(i=1,2,..,m)$ and $y^{k}=\bar{y}^{k}$, then $w^{k+1}$ produced by Algorithm 1 is a solution of $VI(\mathcal{W},\mathcal{F},\vartheta)$. \\
\indent In the following, we show that the accuracy of $\bar{w}^{k}$ to a solution of $VI(\mathcal{W},\mathcal{F},\vartheta)$ is measured by the quantity $(w-\bar{w}^{k})^{\top}\mathcal{Q}(u^{k}-\bar{w}^{k})$.\\
\textbf{Lemma 3.2} Let the sequence $\{w^{k}\}$ be generated by Algorithm. Then we have
\begin{equation}
(w-\bar{w}^{k})^{\top}\mathcal{Q}(w^{k}-\bar{w}^{k})=\frac{1}{2}\left(\|w-w^{k+1}\|_{\mathcal{H}}^{2}-\|w-w^{k}\|_{\mathcal{H}}^{2}+\|w^{k}-\bar{w}^{k}\|_{\mathcal{ N}}^{2}\right),~~~~\forall w\in \mathcal{W},
\end{equation}
where $\mathcal{ N}=\mathcal{Q}+\mathcal{Q}^{\top}-\mathcal{M}^{\top}\mathcal{H}\mathcal{M}$.\\
$Proof$ It follows from the positive definiteness of $\mathcal{H}$ that
\begin{equation*}
\begin{split}
\left(w-\bar{w}^{k}\right)^{\top}\mathcal{H}\left(w^{k}-w^{k+1}\right)=&\frac{1}{2}\left(\|w-w^{k+1}\|_{\mathcal{H}}^{2}-\|w-w^{k}\|_{\mathcal{H}}^{2}\right)\\
&+\frac{1}{2}\left(\|w^{k}-\bar{w}^{k}\|_{\mathcal{H}}^{2}-\|w^{k+1}-\bar{w}^{k}\|_{\mathcal{H}}^{2}\right).
\end{split}
\end{equation*}
Furthermore, by the identity (14) and $\mathcal{Q}=\mathcal{H}\mathcal{M}$, we have
\begin{equation}
\begin{split}
\left(w-\bar{w}^{k}\right)^{\top}\mathcal{Q}\left(w^{k}-\bar{w}^{k}\right)=&\frac{1}{2}\left(\|w-w^{k+1}\|_{\mathcal{H}}^{2}-\|w-w^{k}\|_{\mathcal{H}}^{2}\right)\\
&+\frac{1}{2}\left(\|w^{k}-\bar{w}^{k}\|_{\mathcal{H}}^{2}-\|w^{k+1}-\bar{w}^{k}\|_{\mathcal{H}}^{2}\right).
\end{split}
\end{equation}
And
\begin{equation}
\begin{split}
\|w^{k}-\bar{w}^{k}\|_{\mathcal{H}}^{2}-\|w^{k+1}-\bar{w}^{k}\|_{\mathcal{H}}^{2}&=\|w^{k}-\bar{w}^{k}\|_{\mathcal{H}}^{2}-\|(w^{k}-\bar{w}^{k})-(w^{k}-w^{k+1})\|_{\mathcal{H}}^{2}\\
&=\|w^{k}-\bar{w}^{k}\|_{\mathcal{H}}^{2}-\|(w^{k}-\bar{w}^{k})-\mathcal{M}(w^{k}-\bar{w}^{k})\|_{\mathcal{H}}^{2}\\
&=\left(w^{k}-\bar{w}^{k}\right)^{\top}\left(\mathcal{Q}^{\top}+\mathcal{Q}-\mathcal{M}^{\top}\mathcal{H}\mathcal{M}\right)\left(w^{k}-\bar{w}^{k}\right)\\
&=\|w^{k}-\bar{w}^{k}\|_{\mathcal{N}}^{2},
\end{split}
\end{equation}
where $\mathcal{ N}=\mathcal{Q}+\mathcal{Q}^{\top}-\mathcal{M}^{\top}\mathcal{H}\mathcal{M}$. Thus, the assertion (23) is proved directly by (24) and (25).$~~~~~~~~~~~~~~~~~~~~~~~~~~~~~~~~~~~~~~~~~~~~~~~~~~~~~~~~~~~~~~~~~~~~~~~~~~~~~~~~~~~~~~~~~~~~~~~~~~~~~~~~~~~~~~~~~~\square$\\
\indent In the following, we give a crucial inequality for the global convergence of Algorithm 1 that indicates the sequence $\{w^{k}\}$ generated by Algorithm 1 is Fej$\grave{e}$r monotone with respect to the solution $\mathcal{W}^{*}$.\\
\textbf{Lemma 3.3} Let $\{w^{k}\}$ be the sequence generated by Algorithm 1. Then we have
\begin{equation}
\|w^{k+1}-w^{*}\|_{\mathcal{H}}^{2}\leq\|w^{k}-w^{*}\|_{\mathcal{H}}^{2}-\|w^{k}-\bar{w}^{k}\|_{\mathcal{N}}^{2},~~\forall w^{*}\in\mathcal{W}^{*}.
\end{equation}
$Proof$ It follows from (15) and (23) that we obtain
\begin{equation}
\vartheta(u)-\vartheta(\bar{u}^{k})+(w-\bar{w}^{k})^{\top}\mathcal{F}(\bar{w}^{k})\geq\frac{1}{2}\left(\|w-w^{k+1}\|_{\mathcal{H}}^{2}-\|w-w^{k}\|_{\mathcal{H}}^{2}+\|w^{k}-\bar{w}^{k}\|_{\mathcal{ N}}^{2}\right), ~~\forall w\in\mathcal{W}.
\end{equation}
 Furthermore, by setting $w=w^{*}\in\mathcal{W}^{*}$ in (27) and the monotonicity of $\mathcal{F}(\cdot)$, we get
\begin{equation*}
\begin{split}
0&\geq\vartheta(u^{*})-\vartheta(\bar{u}^{k})+(w^{*}-\bar{w}^{k})^{\top}\mathcal{F}(\bar{w}^{k})\\
&\geq\frac{1}{2}\left(\|w^{*}-w^{k+1}\|_{\mathcal{H}}^{2}-\|w^{*}-w^{k}\|_{\mathcal{H}}^{2}+\|w^{k}-\bar{w}^{k}\|_{\mathcal{ N}}^{2}\right),
\end{split}
\end{equation*}
which implies the assertion (26).$~~~~~~~~~~~~~~~~~~~~~~~~~~~~~~~~~~~~~~~~~~~~~~~~~~~~~~~~~~~~~~~~~~~~~~~~~~~~~~~~~\square$\\
\indent Now, we establish the global convergence of Algorithm 1.\\
\textbf{Theorem 3.1} The sequence $\{w^{k}\}$ generated by Algorithm 1 converges to an optimal solution which belongs to $\mathcal{W}^{*}$ for problem (1) from any starting point.\\
$Proof$ It follows from (26) that
\begin{equation}
\|w^{k}-\bar{w}^{k}\|_{\mathcal{N}}^{2}\leq\|w^{k}-w^{*}\|_{\mathcal{H}}^{2}-\|w^{k+1}-w^{*}\|_{\mathcal{H}}^{2},~~\forall w^{*}\in\mathcal{W}^{*}.
\end{equation}
Summing the above inequality over $k=0,1,..,\infty$, we have
\begin{equation*}
\sum_{k=0}^{\infty}\|w^{k}-\bar{w}^{k}\|_{\mathcal{N}}^{2}\leq\|w^{0}-w^{*}\|_{\mathcal{H}}^{2},
\end{equation*}
which implies that
\begin{equation}
\lim_{k\rightarrow\infty}\|w^{k}-\bar{w}^{k}\|_{\mathcal{N}}^{2}=0.
\end{equation}
It follows from (29) that $\bar{w}^{k}\rightarrow w^{k}$ $(k\rightarrow\infty)$. Moreover, by (15) and (29), we have
\begin{equation}
\lim_{k\rightarrow\infty}\left\{\vartheta(u)-\vartheta(\bar{u}^{k})+(w-\bar{w}^{k})^{\top}\mathcal{F}(\bar{w}^{k})\right\}\geq0,~~~~\forall w\in \mathcal{W},
\end{equation}
From (26), we know that the sequence $\{w^{k}\}$ is bounded. Then, it follows from (29) that the sequence $\{\bar{w}^{k}\}$ is also bounded. Thus, the sequence $\{\bar{w}^{k}\}$ has at least one cluster point. We assume that $w^{\infty}$ is a cluster point of $\{\bar{w}^{k}\}$ and the subsequence $\{\bar{w}^{k_{j}}\}$ converges to $w^{\infty}$ (That is, $\bar{w}^{k_{j}}\rightarrow w^{\infty}$ $(j\rightarrow\infty)$ ). From (30), we get
\begin{equation*}
\vartheta(u)-\vartheta(u^{\infty})+(w-w^{\infty})^{\top}\mathcal{F}(w^{\infty})\geq0,~~~~\forall w\in \mathcal{W},
\end{equation*}
which indicates that $w^{\infty}\in\mathcal{W}^{*}$. It follows from (29) that $\lim_{k\rightarrow\infty}\|w^{k}-\bar{w}^{k}\|_{\mathcal{H}}^{2}=0.$
Therefore, for any $l>0$, there exists an integer $N$ such that
\begin{equation}
\|w^{k_{N}}-\bar{w}^{k_{N}}\|_{\mathcal{H}}<\frac{l}{2},~~\mathrm{and}~~\|\bar{w}^{k_{N}}-w^{\infty}\|_{\mathcal{H}}<\frac{l}{2}.
\end{equation}
Furthermore, It follows from (26) and (31) that for any $k\geq k_{N}$
\begin{equation*}
\|w^{k}-w^{\infty}\|_{\mathcal{H}}\leq\|w^{k_{N}}-w^{\infty}\|\leq\|w^{k_{N}}-\bar{w}^{k_{N}}\|_{\mathcal{H}}+\|\bar{w}^{k_{N}}-w^{\infty}\|_{\mathcal{H}}<l,
\end{equation*}
which shows that $\lim_{k\rightarrow\infty}w^{k}=w^{\infty}$. It follows from $w^{\infty}\in\mathcal{W}^{*}$ that we thus complete the proof.\\
$~~~~~~~~~~~~~~~~~~~~~~~~~~~~~~~~~~~~~~~~~~~~~~~~~~~~~~~~~~~~~~~~~~~~~~~~~~~~~~~~~~~~~~~~~~~~~~~~~~~~~~~~~~~~~~~~~~~~~~~~~~~\square$\\
\section{The worst case $\mathcal{O}(1/t)$ convergence rate }
In this section, we establish the worst-case convergence rate for the Algorithm 1 in both ergodic and nonerodic senses.
\subsection{A worst case $\mathcal{O}(1/t)$ convergence rate in a ergodic sense}
In this subsection, we show a worst $\mathcal{O}(1/t)$ convergence rate of the Algorithm 1 in a ergodic sense.\\
\textbf{Theorem 4.1} Let $\{w^{k}\}$ be the sequence generated by the Algorithm 1. For any integer $t>0$, let $w_{t}=\frac{1}{t+1}\sum_{k=0}^{t}\bar{w}^{k}$, then we have $w_{t}\in\mathcal{W}$ and
\begin{equation}
\vartheta(u_{t})-\vartheta(u)+(w_{t}-w)^{\top}\mathcal{F}(w)\leq\frac{1}{2(t+1)}\|w-w^{0}\|_{\mathcal{H}}^{2},~~\forall w\in\mathcal{W}.
\end{equation}
$Proof$ Obviously, we know that $\bar{w}^{k}\in\mathcal{W}$ for all $k\geq0$. It follows from the convexity of $\mathcal{X}_{1}, \mathcal{X}_{2},.., \mathcal{X}_{m}$ that $w_{t}\in\mathcal{W}$. And because of the monotonicity of $\mathcal{F}(\cdot)$, we have
\begin{equation}
\vartheta(u)-\vartheta(\bar{u}^{k})+(w-\bar{w}^{k})^{\top}\mathcal{F}(w)+\frac{1}{2}\|w^{k}-w\|_{\mathcal{H}}^{2}\geq\frac{1}{2}\|w^{k+1}-w\|_{\mathcal{H}}^{2},~~\forall w\in\mathcal{W}.
\end{equation}
Summing the above inequality over $k=0,1,...,t$, we get
\begin{equation}
(t+1)\vartheta(u)-\sum_{k=0}^{t}\vartheta(\bar{u}^{k})+\left((t+1)w-\sum_{k=0}^{t}\bar{w}^{k}\right)^{\top}\mathcal{F}(w)+\frac{1}{2}\|w-w^{0}\|_{\mathcal{H}}^{2}\geq0,~~\forall w\in\mathcal{W}.
\end{equation}
According to the definition of $w_{t}$ and the convexity of $\vartheta(u)$, (34) can be rewritten as
\begin{equation}
\vartheta(u_{t})-\vartheta(u)+(w_{t}-w)^{\top}\mathcal{F}(w)\leq\frac{1}{2(t+1)}\|w-w^{0}\|_{\mathcal{H}}^{2},~~\forall w\in\mathcal{W}.
\end{equation}
We thus complete the proof.$~~~~~~~~~~~~~~~~~~~~~~~~~~~~~~~~~~~~~~~~~~~~~~~~~~~~~~~~~~~~~~~~~~~~~~~~~~~~~~~~~~~~~~\square$\\
\indent Let an compact set $G\subset\mathcal{W}$ and define
\begin{equation*}
g=\mathrm{sup}\left\{\|w-w^{0}\|_{\mathcal{H}}|w\in G\right\},
\end{equation*}
where $w^{0}=\left(x_{1}^{0},x_{2}^{0},...,x_{m}^{0},y^{0}\right)$ is the initial point. After $t$ iterations of Algorithm 1, we can find a $w_{t}\in\mathcal{W}$ such that
\begin{equation*}
\mathrm{sup}_{w\in G}\left\{\vartheta(u_{t})-\vartheta(u)+(w_{t}-w)^{\top}\mathcal{F}(w)\right\}\leq\frac{g^{2}}{2t},
\end{equation*}
which implies that the Algorithm 1 reaches
\begin{equation*}
\vartheta(u_{t})-\vartheta(u)+(w_{t}-w)^{\top}\mathcal{F}(w)\leq\varepsilon, ~~\forall w\in G,
\end{equation*}
at most $t=\frac{g^{2}}{2\varepsilon}$ iterations. We thus establish the worst-case $\mathcal{O}(1/t)$ convergence rate for the Algorithm 1 in a ergodic sense.
\subsection{a worst case $\mathcal{O}(1/t)$ convergence rate in a nonergodic sense}
In this subsection, we give a worst $\mathcal{O}(1/t)$ convergence rate of the Algorithm 1 in a nonergodic sense.
We first have to mention that the term $\|w^{k}-w^{k+1}\|_{\mathcal{H}}^{2}$ can be used to measure the accuracy of an iteration.\\
\textbf{Lemma 4.1} let $\{w^{k+1}\}$ be generated by Algorithm 1, Then $\bar{w}^{k}$ defined in (13) is a solution to (7) if $\|w^{k}-w^{k+1}\|_{\mathcal{H}}^{2}=0$.\\
$Proof$ By Lemma 2.1, Lemma 3.1 and the equation (14), we have
\begin{equation}
\vartheta(u)-\vartheta(\bar{u}^{k})+(w-\bar{w}^{k})^{\top}\mathcal{F}(\bar{w}^{k})\geq(w-\bar{w}^{k})^{\top}\mathcal{H}(w^{k}-w^{k+1}),~~~~\forall w\in \mathcal{W},
\end{equation}
Clearly, $\mathcal{H}$ is positive definite since $\mathcal{G}_{1}$ is positive definite in Algorithm 1 and $\gamma\in(0,2)$. Then the right-hand side of (36) vanishes if $\|w^{k+1}-w^{k}\|_{\mathcal{H}}^{2}=0$, since we can obtain $\mathcal{H}(w^{k+1}-w^{k})=0$ whenever $\|w^{k+1}-w^{k}\|_{\mathcal{H}}^{2}=0$. The assertion is proved.$~~~~~~~~~~~~~~~~~~~~~~~~~~~~~~~~~~~~~~~~~~~~~~~~~~~~~~~~~~~~~~~~~~~~~~~~~~~~~~~~~~~\square$\\
\indent Now, we are ready to establish a worst-case $\mathcal{O}(1/t)$ convergence rate for Algorithm 1 in a nonergodic sense. First, we show the following two lemmas.\\
\textbf{Lemma 4.2} Let the sequence $\{w^{k}\}$ be generated by Algorithm 1 with $\gamma\in(0,2)$ and the associated $\{\bar{w}^{k}\}$ be defined in (13); the matrix $\mathcal{Q}$ be defined in (10). Then, we have
\begin{equation}
(\bar{w}^{k}-\bar{w}^{k+1})^{\top}\mathcal{Q}\left\{(w^{k}-w^{k+1})-(\bar{w}^{k}-\bar{w}^{k+1})\right\}\geq0.
\end{equation}
$Proof$ Setting $w=\bar{w}^{k+1}$ in (13), we obtain
\begin{equation}
\vartheta(\bar{u}^{k+1})-\vartheta(\bar{u}^{k})+(\bar{w}^{k+1}-\bar{w}^{k})^{\top}\mathcal{F}(\bar{w}^{k})\geq(\bar{w}^{k+1}-\bar{w}^{k})^{\top}\mathcal{Q}(w^{k}-\bar{w}^{k}),~~~~\forall w\in \mathcal{W},
\end{equation}
It is easy to know that (13) is also true for $k:=k+1$, that is
\begin{equation*}
\vartheta(u)-\vartheta(\bar{u}^{k+1})+(w-\bar{w}^{k+1})^{\top}\mathcal{F}(\bar{w}^{k+1})\geq(w-\bar{w}^{k+1})^{\top}\mathcal{Q}(w^{k+1}-\bar{w}^{k+1}),~~~~\forall w\in \mathcal{W},
\end{equation*}
By setting $w=\bar{w}^{k}$ in the above inequality, we have
\begin{equation}
\vartheta(\bar{u}^{k})-\vartheta(\bar{u}^{k+1})+(\bar{w}^{k}-\bar{w}^{k+1})^{\top}\mathcal{F}(\bar{w}^{k+1})\geq(\bar{w}^{k}-\bar{w}^{k+1})^{\top}\mathcal{Q}(w^{k+1}-\bar{w}^{k+1}),~~~~\forall w\in \mathcal{W},
\end{equation}
Adding (38) and (39) and using the monotonicity of $\mathcal{F}$, the assertion (37) is proved immediately.\\$~~~~~~~~~~~~~~~~~~~~~~~~~~~~~~~~~~~~~~~~~~~~~~~~~~~~~~~~~~~~~~~~~~~~~~~~~~~~~~~~~~~~~~~~~~~~~~~~~~~~~~~~~~~~~~~~~~~~~~~~~~~~~~~~~\square$
\textbf{Lemma 4.3} Let the sequence $\{w^{k}\}$ be generated by the Algorithm 1 and $\{\bar{w}^{k}\}$ be defined in (13); the matrices $\mathcal{Q}$, $\mathcal{M}$ and $\mathcal{H}$ be defined in (10), (11) and (12), respectively. Then we have
\begin{equation}
\begin{split}
(&w^{k}-\bar{w}^{k})^{\top}\mathcal{M}^{\top}\mathcal{H}\mathcal{M}\left\{(w^{k}-
\bar{w}^{k})-(w^{k+1}-\bar{w}^{k+1})\right\}\\
&\geq\frac{1}{2}\|(w^{k}-\bar{w}^{k})-(w^{k+1}-\bar{w}^{k+1})\|_{(\mathcal{Q}^{\top}+\mathcal{Q})}^{2}.
\end{split}
\end{equation}
$Proof$ Adding the equation
\begin{equation*}
\begin{split}
&\left\{(w^{k}-w^{k+1})-(\bar{w}^{k}-\bar{w}^{k+1})\right\}^{\top}\mathcal{Q}\left\{(w^{k}-
w^{k+1})-(\bar{w}^{k}-\bar{w}^{k+1})\right\}\\
&=\frac{1}{2}\|(w^{k}-\bar{w}^{k})-(w^{k+1}-\bar{w}^{k+1})\|_{(\mathcal{Q}^{\top}+\mathcal{Q})}^{2}
\end{split}
\end{equation*}
to both sides of (37), we have
\begin{equation}
\begin{split}
(&w^{k}-w^{k+1})^{\top}\mathcal{Q}\left\{(w^{k}-
w^{k+1})-(\bar{w}^{k}-\bar{w}^{k+1})\right\}\\
&\geq\frac{1}{2}\|(w^{k}-\bar{w}^{k})-(w^{k+1}-\bar{w}^{k+1})\|_{(\mathcal{Q}^{\top}+\mathcal{Q})}^{2}
\end{split}
\end{equation}
It follows from $\mathcal{Q}=\mathcal{H}\mathcal{M}$, (14) and (41) that the assertion (40) $\mathrm{holds}.~~~~~~~~~~~~~~~~~~~~~\square$
\textbf{Theorem 4.4} Let $\gamma\in(0,2)$, the sequence $\{w^{k}\}$ be generated by Algorithm 1 and the matrix $\mathcal{H}$ be defined in (12). Then, we have
\begin{equation}
\|w^{k+1}-w^{k+2}\|_{\mathcal{H}}^{2}\leq\|w^{k}-w^{k+1}\|_{\mathcal{H}}^{2}.
\end{equation}
$Proof$ Setting $c=\mathcal{M}(w^{k}-\bar{w}^{k})$ and $d=\mathcal{M}(w^{k+1}-\bar{w}^{k+1})$, we get
\begin{equation*}
\|c\|_{\mathcal{H}}^{2}-\|d\|_{\mathcal{H}}^{2}=2c^{\top}\mathcal{H}(c-d)-\|c-d\|_{\mathcal{H}}^{2}.
\end{equation*}
We thus have
\begin{equation}
\begin{split}
\|\mathcal{M}&(w^{k}-\bar{w}^{k})\|_{\mathcal{H}}^{2}-\left\|\mathcal{M}(w^{k+1}-\bar{w}^{k+1})\right\|_{\mathcal{H}}^{2}\\
=&2(w^{k}-\bar{w}^{k})\mathcal{M}^{\top}\mathcal{H}\mathcal{M}\left\{(w^{k}-\bar{w}^{k})-(w^{k+1}-\bar{w}^{k+1})\right\}\\
&-\left\|\mathcal{M}\left\{(w^{k}-\bar{w}^{k})-(w^{k+1}-\bar{w}^{k+1})\right\}\right\|_{\mathcal{H}}^{2}.
\end{split}
\end{equation}
By (40) and (43), we have
\begin{equation}
\begin{split}
\|\mathcal{M}&(w^{k}-\bar{w}^{k})\|_{\mathcal{H}}^{2}-\left\|\mathcal{M}(w^{k+1}-\bar{w}^{k+1})\right\|_{\mathcal{H}}^{2}\\
\geq&\|(w^{k}-\bar{w}^{k})-(w^{k+1}-\bar{w}^{k+1})\|_{(\mathcal{Q}^{\top}+\mathcal{Q})}^{2}\\
&-\left\|\mathcal{M}[(w^{k}-\bar{w}^{k})-(w^{k+1}-\bar{w}^{k+1})]\right\|_{\mathcal{H}}^{2}\\
=&\|(w^{k}-\bar{w}^{k})-(w^{k+1}-\bar{w}^{k+1})\|_{\mathcal{N}}^{2},\\
\end{split}
\end{equation}
where $\mathcal{N}=\mathcal{Q}^{\top}+\mathcal{Q}-\mathcal{M}^{\top}\mathcal{H}\mathcal{M}$.
It follows from $\mathcal{Q}=\mathcal{H}\mathcal{M}$ that
\begin{equation*}
\begin{split}
\mathcal{N}
=&(\mathcal{Q}^{\top}+\mathcal{Q})-\mathcal{M}^{\top}\mathcal{Q}\\
=&\left(                 
  \begin{array}{cccc}   
    2\mathcal{G}_{1} & 0&0\\  
    0&2(p_{m}+\rho \mathcal{A}_{m}^{\top}\mathcal{A}_{m})&-\gamma\mathcal{A}_{m}^{\top}\\
    0&-\gamma \mathcal{A}_{m}&\frac{2}{\rho}I_{\ell}
  \end{array}
\right)
-\left(                 
  \begin{array}{cccc}   
    \mathcal{G}_{1}& 0 & 0\\  
    0 & p_{m}+2\rho \mathcal{A}_{m}^{\top}\mathcal{A}_{m} & -\gamma\mathcal{A}_{m}^{\top}\\  
    0&-\gamma\mathcal{A}_{m}&\frac{\gamma }{\rho}I_{\ell}\\
  \end{array}
\right)\\
=&\left(                 
  \begin{array}{cccc}   
    \mathcal{G}_{1} &0& 0\\  
   0 & P_{m} & 0\\  
   0&0&\frac{2-\gamma}{\rho}I_{\ell}
  \end{array}
\right).
\end{split}
\end{equation*}
By the positive definiteness of $P_{m}$ and $\mathcal{G}_{1}$ , it is easy to see that $\mathcal{N}$ is positive definite when $\gamma\in(0,2)$. Then, we have
\begin{equation*}
\left\|\mathcal{M}(w^{k}-\bar{w}^{k})\right\|_{\mathcal{H}}^{2}-\left\|\mathcal{M}(w^{k+1}-\bar{w}^{k+1})\right\|_{\mathcal{H}}^{2}\geq0.
\end{equation*}
By $w^{k}-w^{k+1}=\mathcal{M}(w^{k}-\bar{w}^{k})$, the assertion (42) follows $\mathrm{immediately}.~~~~~~~~~~~~~~\square$\\
\indent In order to further finish the proof of a worst-case $\mathcal{O}(1/t)$ convergence rate for Algorithm 1 in a nonergodic sense, we give the following two lemmas.\\
\textbf{Lemma 4.4} Let $\{x_{m}^{k}\}$ be the sequence generated by Algorithm 1 with $\gamma\in(0,2)$. Then, we have
\begin{equation}
(x_{m}^{k}-x_{m}^{k+1})^{\top}\mathcal{A}_{m}^{\top}(y^{k}-y^{k+1})\geq\frac{1}{2}\|x_{m}^{k}-x_{m}^{k+1}\|_{P_{m}}^{2}-\frac{1}{2}\|x_{m}^{k-1}-x_{m}^{k}\|_{P_{m}}^{2}.
\end{equation}
$Proof$ By the optimality condition of the $x_{m}$-subproblem in lemma 3.1 and (6), there exists $\xi\in\partial \vartheta_{m}(x_{m}^{k+1})$ such that
\begin{equation}
(x_{m}-x_{m}^{k+1})^{\top}\left\{\xi-\mathcal{A}_{m}^{\top}y^{k+1}+P_{m}(x_{m}^{k+1}-x_{m}^{k})\right\}\geq0,~~~~\forall~x_{m}\in\mathcal{X}_{m},
\end{equation}
where $\partial \vartheta_{m}(x_{m})$ is a subdifferential of $\vartheta_{m}(x_{m})$. Setting $x_{m}=x_{m}^{k}$ in (46), then we have
\begin{equation}
(x_{m}^{k}-x_{m}^{k+1})^{\top}\left\{\xi-\mathcal{A}_{m}^{\top}y^{k+1}+P_{m}(x_{m}^{k+1}-x_{m}^{k})\right\}\geq0.
\end{equation}
Furthermore, setting $k:=k+1$ in (46), there exists $\zeta\in\partial \vartheta_{m}(x_{m}^{k})$ such that
\begin{equation}
(x_{m}-x_{m}^{k})^{\top}\left\{\zeta-\mathcal{A}_{m}^{\top}y^{k}+P_{m}(x_{m}^{k}-x_{m}^{k-1})\right\}\geq0.
\end{equation}
Similarly, by $x_{m}=x_{m}^{k+1}$ in (48), we have
\begin{equation}
(x_{m}^{k+1}-x_{m}^{k})^{\top}\left\{\zeta-\mathcal{A}_{m}^{\top}y^{k}+P_{m}(x_{m}^{k}-x_{m}^{k-1})\right\}\geq0.
\end{equation}
Adding (47) to (49) and using the monotonicity of the $\partial \vartheta_{m}(\cdot)$, we get
\begin{equation*}
\begin{split}
(x_{m}^{k}-x_{m}^{k+1})^{\top}\mathcal{A}_{m}^{\top}(y^{k}-y^{k+1})&\geq(x_{m}^{k+1}-x_{m}^{k})^{\top}P_{m}(x_{m}^{k+1}-x_{m}^{k}+x_{m}^{k-1}-x_{m}^{k})\\
&=\|x_{m}^{k+1}-x_{m}^{k}\|_{P_{m}}^{2}+(x_{m}^{k+1}-x_{m}^{k})^{\top}P_{m}(x_{m}^{k-1}-x_{m}^{k})
\end{split}
\end{equation*}
By the inequality
\begin{equation*}
 (x_{m}^{k+1}-x_{m}^{k})^{\top}P_{m}(x_{m}^{k-1}-x_{m}^{k})\geq-\frac{1}{2}\|x_{m}^{k}-x_{m}^{k+1}\|_{P_{m}}^{2}-\frac{1}{2}\|x_{m}^{k-1}-x_{m}^{k}\|_{P_{m}}^{2}, \end{equation*}
the assertion (45) is proved.$~~~~~~~~~~~~~~~~~~~~~~~~~~~~~~~~~~~~~~~~~~~~~~~~~~~~~~~~~~~~~~~~~~~~~~~~~~~~~~~~~~~\square$\\
\textbf{Lemma 4.5} Let the sequence $\{w^{k}\}$ be generated by the Algorithm 1 with $\gamma\in(0,2)$ and the associated $\{\bar{w}^{k}\}$ be defined in (13), then there exists $0<\sigma_{\gamma}\leq1$ such that
\begin{equation}
\begin{split}
\|&w^{k}-\bar{w}^{k}\|_{\mathcal{N}}^{2}\\
&\geq\sigma_{\gamma}\left(\|R^{k}-\bar{R}^{k}\|_{\mathcal{G}_{1}}^{2}+\|x_{m}^{k}-\bar{x}_{m}^{k}\|_{P_{m}}^{2}+\frac{\gamma}{\rho}\|y^{k}-\bar{y}^{k}\|^{2}\right),
\end{split}
\end{equation}
where
\begin{equation}       
\mathcal{N}=\mathcal{Q}^{\top}+\mathcal{Q}-\mathcal{M}^{\top}\mathcal{H}\mathcal{M},~~~~
R^{k}=\left(                 
  \begin{array}{ccc}   
    x_{1}^{k}\\  
     x_{2}^{k}\\  
     \vdots\\
      x_{m-1}^{k}\\
  \end{array}
\right)                 
\end{equation}
$Proof.$ By the definition of $\mathcal{Q}, \mathcal{M}$, $\mathcal{H}$, $\mathcal{G}_{1}$ and $R^{k}$, we have
\begin{equation}
\begin{split}
\|w&^{k}-\bar{w}^{k}\|_{\mathcal{N}}^{2}\\
=&\|R^{k}-\bar{R}^{k}\|_{\mathcal{G}_{1}}^{2}+\|x_{m}^{k}-\bar{x}_{m}^{k}\|_{P_{m}}^{2}+\frac{2-\gamma}{\rho}\|y^{k}-\bar{y}^{k}\|^{2}\\
\geq&\mathrm{min}\left\{\frac{2-\gamma}{\gamma},1\right\}\left(\|R^{k}-\bar{R}^{k}\|_{\mathcal{G}_{1}}^{2}+\|x_{m}^{k}-\bar{x}_{m}^{k}\|_{P_{m}}^{2}+\frac{\gamma}{\rho}\|y^{k}-\bar{y}^{k}\|^{2}\right).
\end{split}
\end{equation}
Let $\sigma_{\gamma}=\mathrm{min}\left\{\frac{2-\gamma}{\gamma},1\right\}$, then it follows from (52) that the assertion (50) $\mathrm{holds}. \square$\\
\indent Finally, we establish a worst-case $\mathcal{O}(1/t)$ convergence rate for the Algorithm 1 in a nonergodic sense.\\
\textbf{Theorem 4.5} Let the sequence $\{w^{k}\}$ be generated by Algorithm 1 with $\gamma\in(0,2)$. Then we have
\begin{equation}
\|w^{t}-w^{t+1}\|_{\mathcal{H}}^{2}\leq\frac{1}{t}\left(\frac{1}{\sigma_{\gamma}}\|w^{0}-w^{*}\|_{\mathcal{H}}^{2}+\|x_{m}^{0}-x_{m}^{1}\|_{\mathcal{H}}^{2}\right).
\end{equation}
$Proof$ By $w=w^{*}$ in (28), we have
\begin{equation}
\begin{split}
\vartheta(u^{*})&-\vartheta(\bar{u}^{k})+(w^{*}-\bar{w}^{k})^{\top}\mathcal{F}(w^{*})\\
\geq& \frac{1}{2}(\|w^{*}-w^{k+1}\|_{\mathcal{H}}^{2}-\|w^{*}-w^{k}\|_{\mathcal{H}}^{2})+\frac{1}{2}\|w^{k}-\bar{w}^{k}\|_{\mathcal{N}}^{2},
\end{split}
\end{equation}
where $\mathcal{N}=\mathcal{Q}^{\top}+\mathcal{Q}-\mathcal{M}^{\top}\mathcal{H}\mathcal{M}$.
By (7) and (54), we can obtain
\begin{equation*}
\|w^{k}-\bar{w}^{k}\|_{\mathcal{N}}^{2}\leq\|w^{*}-w^{k}\|_{\mathcal{H}}^{2}-\|w^{*}-w^{k+1}\|_{\mathcal{H}}^{2}.
\end{equation*}
According to the Theorem 4.1, By the positive definiteness of $P_{m}$ and $\mathcal{G}_{1}$, it is easy to see that $\mathcal{N}$ is positive definite. Thus, we have
\begin{equation}
\sum_{k=0}^{\infty}\|w^{k}-\bar{w}^{k}\|_{\mathcal{N}}^{2}\leq\|w^{0}-w^{*}\|_{\mathcal{H}}^{2}.
\end{equation}
Moreover, by the definition of $\mathcal{H}$ in (12), we have
\begin{equation*}
\begin{split}
&\|w^{k}-w^{k+1}\|_{\mathcal{H}}^{2}\\
=&\|R^{k}-\bar{R}^{k}\|_{\mathcal{G}_{1}}^{2}+\|x_{m}^{k}-\bar{x}_{m}^{k}\|_{P_{m}}^{2}\\
&+\frac{1}{\rho\gamma}(\|\rho \mathcal{A}_{m}(x_{m}^{k}-x_{m}^{k+1})\|^{2}+\|y^{k}-y^{k+1}\|^{2}
+2(1-\gamma)\rho(x_{m}^{k}-x_{m}^{k+1})^{\top}\mathcal{A}_{m}^{\top}(y^{k}-y^{k+1}))\\
=&\|R^{k}-\bar{R}^{k}\|_{\mathcal{G}_{1}}^{2}+\|x_{m}^{k}-\bar{x}_{m}^{k}\|_{P_{m}}^{2}+\frac{1}{\rho\gamma}\|\rho \mathcal{A}_{m}(x_{m}^{k}-x_{m}^{k+1})+(y^{k}-y^{k+1})\|^{2}\\
&-2(x_{m}^{k}-x_{m}^{k+1})^{\top}\mathcal{A}_{m}^{\top}(y^{k}-y^{k+1}).
\end{split}
\end{equation*}
Furthermore, by  (6) and (13), we have
\begin{equation}
y^{k}-y^{k+1}=-\rho \mathcal{A}_{m}(x_{m}^{k}-\bar{x}_{m}^{k})+\gamma(y^{k}-\bar{y}^{k}).
\end{equation}
We thus get
\begin{equation}
\begin{split}
\|w^{k}-w^{k+1}\|_{\mathcal{H}}^{2}
=&\|R^{k}-\bar{R}^{k}\|_{\mathcal{G}_{1}}^{2}+\|x_{m}^{k}-\bar{x}_{m}^{k}\|_{P_{m}}^{2}\\
&+\frac{\gamma}{\rho}\|y^{k}-\bar{y}^{k}\|^{2}-2(x_{m}^{k}-x_{m}^{k+1})^{\top}\mathcal{A}_{m}^{\top}(y^{k}-y^{k+1}).
\end{split}
\end{equation}
It follows from (45), (50), (55) and (56) that
\begin{equation}
\begin{split}
\sum_{k=1}^{t}\|w^{k}-w^{k+1}\|_{\mathcal{H}}^{2}\leq&\frac{1}{\sigma_{\gamma}}\sum_{k=1}^{t}\|w^{k}-\bar{w}^{k}\|_{\mathcal{N}}^{2}\\
&+\sum_{k=1}^{t}\left(\|x_{m}^{k-1}-x_{m}^{k}\|_{P_{m}}^{2}-\|x_{m}^{k}-x_{3}^{k+1}\|_{P_{m}}^{2}\right)\\
\leq&\frac{1}{\sigma_{\gamma}}\|w^{0}-w^{*}\|_{\mathcal{H}}^{2}+\|x_{m}^{0}-x_{m}^{1}\|_{P_{m}}^{2}.
\end{split}
\end{equation}
According to the Theorem 4.1, the sequence $\{\|w^{k}-w^{k+1}\|_{\mathcal{H}}^{2}\}$ is non-increasing. Thus, we get
\begin{equation}
\begin{split}
t\|w^{t}-w^{t+1}\|_{\mathcal{H}}^{2}&\leq\sum_{k=1}^{t}\|w^{k}-w^{k+1}\|_{\mathcal{H}}^{2}\\
&\leq\frac{1}{\sigma_{\gamma}}\|w^{0}-w^{*}\|_{\mathcal{H}}^{2}+\|x_{m}^{0}-x_{m}^{1}\|_{\mathcal{H}}^{2}.
\end{split}
\end{equation}
 Based on the above inequality, the assertion (53) is proved.$~~~~~~~~~~~~~~~~~~~~~~~~~~~~~~~~~~~~~~~~~\square$\\
\\
\indent It follows from the conclusion in Theorem 4.5 that we can obtain a worst-case $\mathcal{O}(1/t)$ convergence rate for Algorithm 1 in a nonergodic sense. .
\section{Numerical experiments}
In this subsection, we report some numerical result for the proposed method by calibrating the correlation matrices. We also compare its numerical performance with the ADMM with Gaussian back substitution (ADMM-G) [20] and an algorithm twisted from generalized ADMM (TADMM) [32].\\
\indent All experiments are implemented in MATLAB R2010b on a hp-notebook with an Intel Core i5-3340M CPU at 2.70 GHz and 8 GB memory.\\
\indent We first consider to solve the following matrix optimization problem
\begin{equation}
\mathrm{min}~\{\frac{1}{2}\|X-C\|_{F}^{2}|X\in S_{+}^{n}\cap \mathcal{S}_{B}\},
\end{equation}
where
\begin{equation*}
S_{+}^{n}=\{H\in R^{n\times n}|H^{\top}=H,H\succeq\textbf{0}\},
\end{equation*}
and
\begin{equation*}
S_{B}=\{H\in R^{n\times n}|H^{\top}=H,H_{L}\preceq H\preceq H_{U}\}.
\end{equation*}
We set  $\mathcal{A}_{1}= \left(\begin{array}{ccc}   
    I\\  
    I\\  
    \textbf{0}
  \end{array}\right)$, $\mathcal{A}_{2}= \left(\begin{array}{ccc}   
    -I\\  
    \textbf{0}\\
    I
  \end{array}\right)$, $\mathcal{A}_{3}= \left(\begin{array}{ccc}   
      \textbf{0}\\
    -I\\  
    -I
  \end{array}\right)$. Then, the above problem (60) could be converted to the following equivalent form
\begin{equation}
\begin{array}{l}
\mathrm{min}~\frac{1}{2}\|X_{1}-C\|^{2}+\frac{1}{2}\|X_{2}-C\|^{2}+\frac{1}{2}\|X_{3}-C\|^{2}\\
\mathrm{s.t.}   ~~~~\mathcal{A}_{1}X_{1}+\mathcal{A}_{2}X_{2}+\mathcal{A}_{3}X_{3}=\textbf{0},\\
~~~~~~~~X_{1},~X_{2}\in S_{+}^{n},~ X_{3}\in S_{B}.
\end{array}
\end{equation}
\indent We then derive the subproblems when Algorithm 1 is applied to (61):
\begin{equation}
\begin{split}
\left\{
\begin{aligned}
X_{1}^{k+1}=&\mathrm{argmin}\{\frac{1}{2}\|X_{1}-C\|^{2}-X_{1}^{\top}\mathcal{A}_{1}^{\top}y^{k}+\frac{\rho}{2}\| \mathcal{A}_{1}X_{1}+\mathcal{A}_{2}X_{2}^{k}+\mathcal{A}_{3}X_{3}^{k}-b \|^{2}\\
&+\frac{1}{2}\|X_{1}-X_{1}^{k}\|_{P_{1}}^{2}|X_{1}\in\mathcal{X}_{1}\} , \\
X_{2}^{k+1}=&\mathrm{argmin}\{\frac{1}{2}\|X_{2}-C\|^{2}-X_{2}^{\top}\mathcal{A}_{2}^{\top}y^{k}+\frac{\rho}{2}\| \mathcal{A}_{1}X_{1}+\mathcal{A}_{2}X_{2}^{k}+\mathcal{A}_{3}X_{3}^{k}-b \|^{2}\\
&+\frac{1}{2}\|X_{2}-X_{2}^{k}\|_{P_{2}}^{2}|X_{2}\in\mathcal{X}_{2}\} , \\
X_{3}^{k+1}=&\mathrm{argmin}\{\frac{1}{2}\|X_{3}-C\|^{2}-X_{3}^{\top}\mathcal{A}_{3}^{\top}y^{k}+\frac{\rho}{2}\| \gamma (\mathcal{A}_{1}X_{1}^{k+1}+\mathcal{A}_{2}X_{2}^{k+1})+(1-\gamma)(b-\mathcal{A}_{3}X_{3}^{k})\\
&+\mathcal{A}_{3}X_{3}-b \|^{2}+\frac{1}{2}\|X_{3}-X_{3}^{k}\|_{P_{3}}^{2}|X_{3}\in\mathcal{X}_{3}\} , \\
y^{k+1}=&y^{k}-\rho\left( \gamma (\mathcal{A}_{1}X_{1}^{k+1}+\mathcal{A}_{2}X_{2}^{k+1})+(1-\gamma)(b-\mathcal{A}_{3}X_{3}^{k})+\mathcal{A}_{3}X_{3}^{k+1}-b\right).\\
\end{aligned}
\right.
\end{split}
\end{equation}
For the $X_{1}$-subproblem in (62), by first-order optimal condition, we have
\begin{equation}
X_{1}^{k+1}=P_{\mathcal{S}_{+}^{n}}\left\{(P_{1}+(1+2\beta)\cdot I)^{-1}\left(-\beta \mathcal{A}_{1}^{\top}\mathcal{A}_{2}X_{2}^{k}-\beta \mathcal{A}_{1}^{\top}\mathcal{A}_{3}X_{3}^{k}+\mathcal{A}_{1}^{\top}y^{k}+C+P_{1}X_{1}^{k}\right)\right\}.
\end{equation}
where $P_{\mathcal{S}_{+}^{n}}(\mathcal{A})=\mathcal{U}\Lambda^{+}\mathcal{U}^{\top}$ and $\Lambda^{+}=\mathrm{max}(\Lambda,0)$, $[\mathcal{U},\Lambda]=\mathrm{eig}(\mathcal{A})$. For the $X_{2}$-subproblem in (62), by simple calculation, we get its solution is
\begin{equation}
X_{2}^{k+1}=P_{\mathcal{S}_{+}^{n}}\left\{(P_{2}+(1+2\beta)\cdot I)^{-1}\left(-\beta \mathcal{A}_{2}^{\top}\mathcal{A}_{1}X_{1}^{k+1}-\beta \mathcal{A}_{2}^{\top}\mathcal{A}_{3}X_{3}^{k}+\mathcal{A}_{2}^{\top}y^{k}+C+P_{2}X_{2}^{k}\right)\right\}.
\end{equation}
The $X_{3}$-subproblem can be solved explicitly via
\begin{equation}
\begin{split}
X_{3}^{k+1}=&P_{\mathcal{S}_{B}}\{(P_{3}+(1+2\beta)\cdot I)^{-1}
(-\beta\gamma(\mathcal{A}_{3}^{\top}\mathcal{A}_{1}*X_{1}^{k+1}+\mathcal{A}_{3}^{\top}\mathcal{A}_{2}X_{2}^{k+1})\\
&+\beta(1-\gamma)\mathcal{A}_{3}^{\top}\mathcal{A}_{3}X_{3}^{k}+\mathcal{A}_{3}^{\top}y^{k}+C+P_{3}X_{3}^{k})\}
\end{split}
\end{equation}
where $S_{B}=\{H\in R^{n\times n}|H_{L}\preceq H\preceq H_{U}\}$ and $P_{\mathcal{S}_{B}}(\mathcal{A})=\mathrm{min}(\mathrm{max}(H_{L},\mathcal{A}),H_{U})$.\\
\indent To implement Algorithm 1, we use the stopping criterion
\begin{equation*}
\mathrm{max}\left\{\frac{\|X_{1}^{k+1}-X_{1}^{k}\|}{\|X_{1}^{1}-X_{1}^{0}\|}, \frac{\|X_{2}^{k+1}-X_{2}^{k}\|}{\|X_{2}^{1}-X_{2}^{0}\|}, \frac{\|X_{3}^{k+1}-X_{3}^{k}\|}{\|X_{3}^{1}-X_{3}^{0}\|}, \frac{\|y^{k+1}-y^{k}\|}{\|y^{1}-y^{0}\|}\right\}<\mathrm{tol},
 \end{equation*}
where $\mathrm{tol}=10^{-6}$. We set $P_{1}=P_{2}=P_{3}=\frac{1}{2}I$, $C=\mathrm{rand}(n,n)$, $C=(C'+C)-\mathrm{ones}(n,n) + eye(n)$, $HU=\mathrm{ones}(n,n)\cdot0.1$ and $HL=-HU$.\\
\indent We first test the sensitivity of $\gamma$ for the L-GADMM (6) and take the matrix optimization model (60) with $n=50$. And we fix $\beta=1$ and choose different values of $\gamma$ in the interval $(0,2)$. We repeat each scheme ten times and report some numerical results for the above experiment, and we plot them in Fig. 5.1.\\
\begin{figure}[H]
  \centering
  \includegraphics[width=1\textwidth]{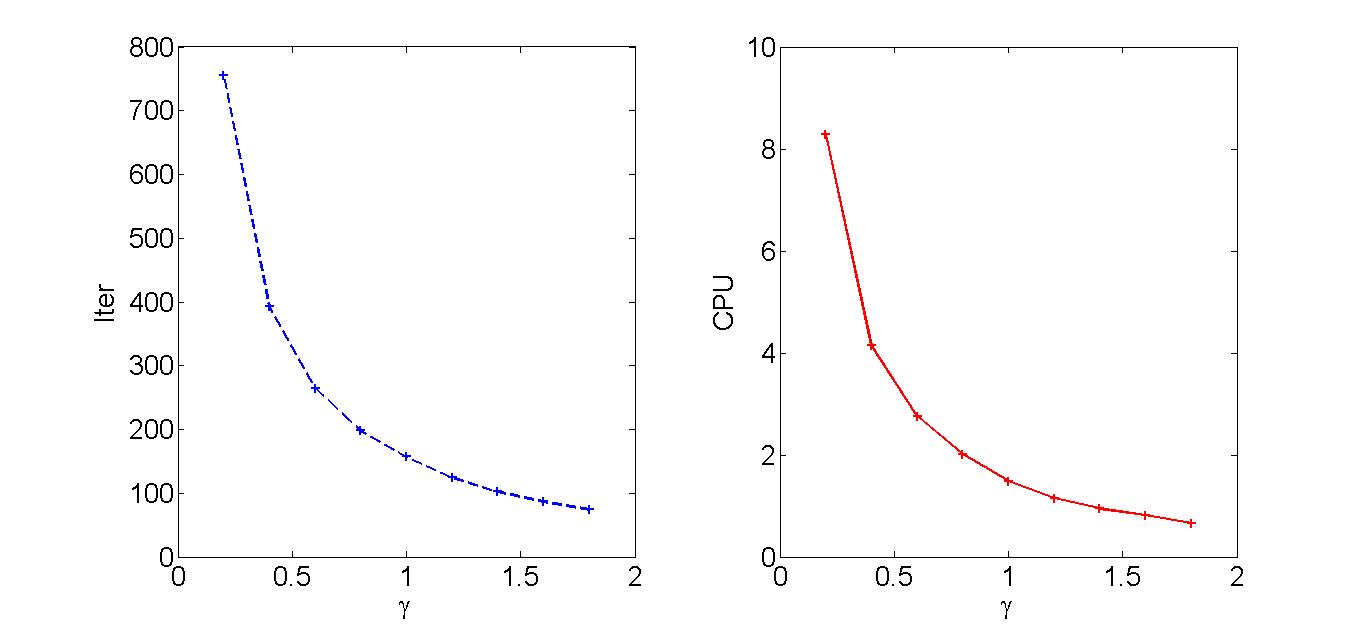}\\
  {Fig.1 Sensitivity test on the relaxation factor $\gamma$}
\end{figure}
From Fig.1, the bigger $\gamma$ often lead to better numerical results since we can see that both the number of iterations and computing time in seconds have decreasing tendency as the relaxation factor $\gamma$ increases. Then, we compare Algorithm 1 with ADMM-G [20], TADMM [32] for their averaged performances. To further observe the convergence of the three tested algorithms, in Fig. 5.2, we plot the evolutions of the objective function value when the three methods are applied to solve (61).
\begin{figure}[H]
  \centering
  \includegraphics[width=1\textwidth]{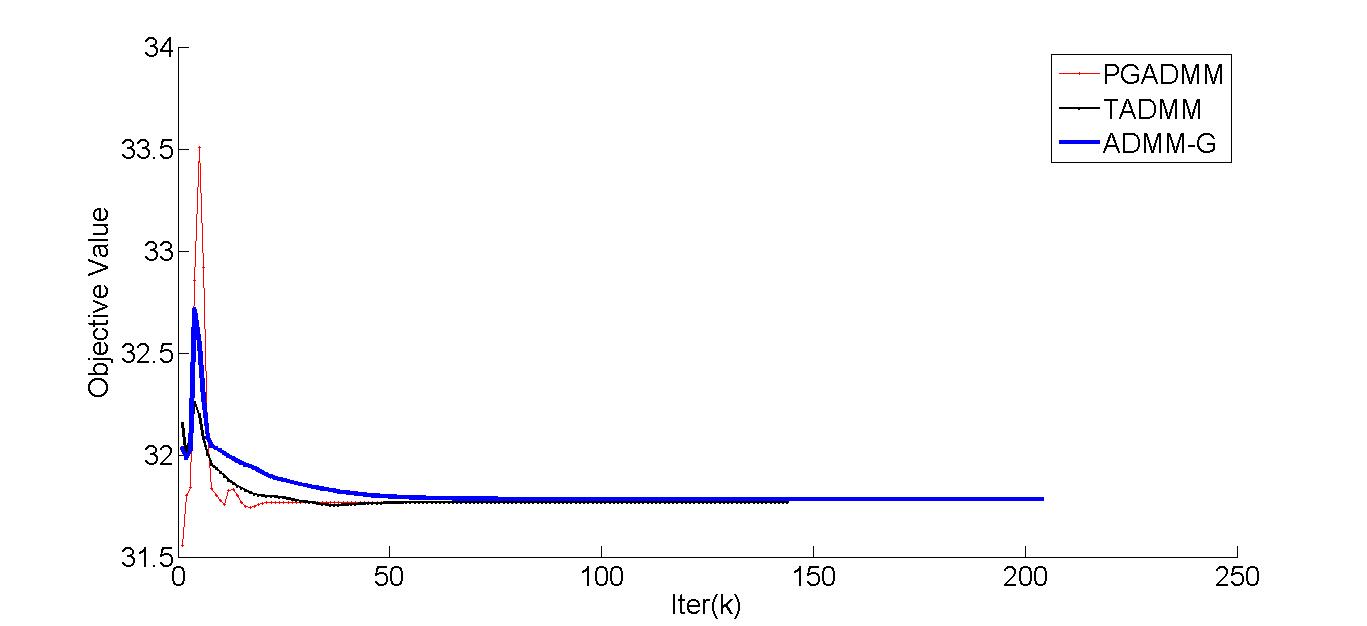}\\
  {Fig.2 Evolutions of the objective function value and the number of iterations}
\end{figure}
From Fig.2, it is easy to see that the proposed algorithm provides faster convergence than that provided by ADMM-G and TADMM; and the figure also shows that the Algorithm 1 has better objective function than that provided by ADMM-G.\\
\begin{table}[H]
\centering  
\begin{tabular}{lccccccccccc}  
\hline
&~~ Alg 1 ($\gamma=1$) &~~Alg 1 ($\gamma=1.9$)&~~TADMM&~~ADMM-G \\ \hline  
n=100\\         
~~Iterations&234&110&217&905\\
~~CPU times(s)&8.302651&4.198581& 7.945454&33.137704\\        
~~Objective value&71.43&71.43&71.43&71.69\\
~~Epsilon&9.988e-007&7.122e-007&9.915e-007&9.989e-007\\
n=300\\
~~Iterations&369&196&345&448\\
~~CPU times(s)&179.435412&100.498216&158.495029&213.224576\\        
~~Objective value&256.90&256.90&256.90& 257.17\\
~~Epsilon&8.148e-007 &9.681e-007&9.992e-007&9.940e-007\\
n=500\\
~~Iterations&428&222&395&505\\
~~CPU times(s)&856.182117&435.490322&732.217027&941.7561791\\        
~~Objective value&430.70&430.70&430.70&430.91\\
~~Epsilon&9.958e-007&9.821e-007 &9.957e-007&9.998e-007
\\\hline
\end{tabular}
\caption{ Numerical comparison between the averaged performance of} Algorithm 1 and TADMM, ADMM-G\label{fig:1}
\end{table}
\indent Furthermore, we give some numerical results under different dimensions (i.e. n=100, 300, 500) in Table 1. The data in this table shows that achieving the same level of objective function values, Algorithm 1 with $\gamma=1.9$ is faster than Algorithm 1 with $\gamma=1$; and both the Algorithm 1 with $\gamma=1.9$ and Algorithm 1 with $\gamma=1$ are much faster than ADMM-G. Clearly, our proposed method with $\gamma=1.9$ are also much faster than TADMM. Moreover, we also give the iterations, computing time in seconds and the epsilon when the stopping criterion is achieved.

According to the above numerical experiment, therefore, the Algorithm 1 is an efficient method for the matrix optimization problem (61).

\section{Conclusions}
In this paper, we propose a linearized GADMM to solve multi-block separable convex optimization problem. Under some mild conditions, we give the proof of the global convergence of Algorithm 1 and establish the worst-case $\mathcal{O}(1/t)$ convergence rate in both the ergodic and nonergodic senses. Our numerical experiment shows that the proposed method is efficient for solving the matrix optimization problem.



\vskip 6mm
\noindent{\bf Appendices: manuscript data}\\
~
\noindent{\bf Appendix A: Comparison of three algorithms in different dimensions}

1. n=100\\

(1)PGADMM10 scheme    n= 100\\
PGADMM10 scheme   k=  234     eps=9.988e-007\\
Elapsed time is 8.302651 seconds.\\

(2)PGADMM11 scheme    n= 100\\
PGADMM11 scheme   k=  110     eps=7.122e-007\\
Elapsed time is 4.198581 seconds.\\

(3)TADMM   n= 100\\
TADMM  k=  217     eps=9.915e-007\\
Elapsed time is 7.945454 seconds.\\

(4)ADMM with Gaussian back substitution   n= 100\\
ADMM with Gaussian back substitution  k=  905     eps=9.989e-007\\
Elapsed time is 33.137704 seconds.\\

2.n=300\\

(1)PGADMM10 scheme    n= 300\\
PGADMM10 scheme   k=  369     eps=8.148e-007\\
Elapsed time is 179.435412 seconds.\\

(2)PGADMM11 scheme    n= 300\\
PGADMM11 scheme   k=  196     eps=9.681e-007\\
Elapsed time is 100.498216 seconds.\\

(3)TADMM   n= 300\\
TADMM  k=  345     eps=9.992e-007\\
Elapsed time is 158.495029 seconds.\\

(4)ADMM with Gaussian back substitution   n= 300\\
ADMM with Gaussian back substitution  k=  448     eps=9.940e-007\\
Elapsed time is 213.224576 seconds.\\

3.n=500\\

(1)PGADMM10 scheme    n= 500\\
PGADMM10 scheme   k=  428     eps=9.958e-007\\
Elapsed time is 856.182117 seconds.\\

(2)PGADMM11 scheme    n= 500\\
PGADMM11 scheme   k=  222     eps=9.821e-007\\
Elapsed time is 435.490322 seconds.\\

(3)TADMM   n= 500\\
TADMM  k=  395     eps=9.957e-007\\
Elapsed time is 732.217027 seconds.\\

(4)ADMM with Gaussian back substitution   n= 500\\
ADMM with Gaussian back substitution  k=  505     eps=9.998e-007\\
Elapsed time is 941.756179 seconds.\\

\noindent{\bf Appendix B: Relationship between CPU time, number of iterations, and parameter Gamma}

1.PGADMM1 scheme    n=  50\\
PGADMM2 scheme   k=  755     eps=9.999e-007\\
Elapsed time is 8.285959 seconds.\\

2.PGADMM2 scheme    n=  50\\
PGADMM2 scheme   k=  393     eps=9.999e-007\\
Elapsed time is 4.141057 seconds.\\

3.PGADMM3 scheme    n=  50\\
PGADMM3 scheme   k=  265     eps=9.958e-007\\
Elapsed time is 2.756718 seconds.\\

4.PGADMM scheme    n=  50\\
PGADMM4 scheme   k=  198     eps=9.878e-007\\
Elapsed time is 2.011761 seconds.\\

5.PGADMM5 scheme    n=  50\\
PGADMM5 scheme   k=  157     eps=9.794e-007\\
Elapsed time is 1.481666 seconds.\\

6.PGADMM6 scheme    n=  50\\
PGADMM6 scheme   k=  124     eps=9.966e-007\\
Elapsed time is 1.152508 seconds.\\

7.PGADMM7 scheme    n=  50\\
PGADMM7 scheme   k=  102     eps=9.688e-007\\
Elapsed time is 0.950794 seconds.\\

8.PGADMM8 scheme    n=  50\\
PGADMM8 scheme   k=   87     eps=8.249e-007\\
Elapsed time is 0.822915 seconds.\\

9.PGADMM9 scheme    n=  50\\
PGADMM9 scheme   k=   74     eps=8.863e-007\\
Elapsed time is 0.658812 seconds.\\

\end{document}